\documentstyle[12pt]{article}
\begin{document}
\newpage
\pagestyle{empty}
\setcounter{page}{0}
%
%%% ** start of amsfont definitions **
\newfont{\twelvemsb}{msbm10 scaled\magstep1}
\newfont{\eightmsb}{msbm8} \newfont{\sixmsb}{msbm6} \newfam\msbfam
\textfont\msbfam=\twelvemsb \scriptfont\msbfam=\eightmsb
\scriptscriptfont\msbfam=\sixmsb \catcode`\@=11
\def\Bbb{\ifmmode\let\next\Bbb@\else \def\next{\errmessage{Use
      \string\Bbb\space only in math mode}}\fi\next}
\def\Bbb@#1{{\Bbb@@{#1}}} \def\Bbb@@#1{\fam\msbfam#1}
\newfont{\twelvegoth}{eufm10 scaled\magstep1}
\newfont{\tengoth}{eufm10} \newfont{\eightgoth}{eufm8}
\newfont{\sixgoth}{eufm6} \newfam\gothfam
\textfont\gothfam=\twelvegoth \scriptfont\gothfam=\eightgoth
\scriptscriptfont\gothfam=\sixgoth \def\frak{\frak@}
\def\frak@#1{{\fam\gothfam{{#1}}}} \def\frak@@#1{\fam\gothfam#1}
\catcode`@=12
%%% ** end of amsfont definitions **
%
%
%
%\def\Bbb{\bf}
\def\CC{{\Bbb C}}
\def\NN{{\Bbb N}}
\def\QQ{{\Bbb Q}}
\def\RR{{\Bbb R}}
\def\ZZ{{\Bbb Z}}
\def\cA{{\cal A}}          \def\cB{{\cal B}}          \def\cC{{\cal C}}
\def\cD{{\cal D}}          \def\cE{{\cal E}}          \def\cF{{\cal F}}
\def\cG{{\cal G}}          \def\cH{{\cal H}}          \def\cI{{\cal I}}
\def\cJ{{\cal J}}          \def\cK{{\cal K}}          \def\cL{{\cal L}} 
\def\cM{{\cal M}}          \def\cN{{\cal N}}          \def\cO{{\cal O}}
\def\cP{{\cal P}}          \def\cQ{{\cal Q}}          \def\cR{{\cal R}} 
\def\cS{{\cal S}}          \def\cT{{\cal T}}          \def\cU{{\cal U}}
\def\cV{{\cal V}}          \def\cW{{\cal W}}          \def\cX{{\cal X}}
\def\cY{{\cal Y}}          \def\cZ{{\cal Z}}
\def\qed{\hfill \rule{5pt}{5pt}}
\newtheorem{lemma}{Lemma}
\newtheorem{prop}{Proposition}
\newtheorem{theo}{Theorem}
\newenvironment{result}{\vspace{.2cm} \em}{\vspace{.2cm}}

$$
\;
$$
\rightline{CPTH-S455.0896}

\vfill
\vfill
\begin{center}

  {\LARGE {\bf {\sf Irreducible Representations of Jordanian 
Quantum  Algebra ${\cal U}_h(sl(2))$ Via a Nonlinear Map}}}
 \\[2cm]

\smallskip 

{\large B. Abdesselam\footnote{abdess@orphee.polytechnique.fr},
A. Chakrabarti\footnote{chakra@orphee.polytechnique.fr} and 
R. Chakrabarti\footnote{Permanent address: Department of Theoretical Physics, University of Madras, Guindy Campus, Madras-600025, India}}

\smallskip 

\smallskip 

\smallskip 

{\em  \footnote{Laboratoire Propre du CNRS UPR A.0014}Centre de Physique 
Th\'eorique, Ecole Polytechnique, \\
91128 Palaiseau Cedex, France.}

\end{center}

\vfill

\begin{abstract}
\noindent The generators of the Jordanian quantum algebra ${\cal U}_h(sl(2))$
are expressed as nonlinear invertible functions of the classical $sl(2)$
generators. This permits immediate explicit construction of the finite dimensional 
irreducible representations of the algebra ${\cal U}_h(sl(2))$.  Using this 
construction, new finite dimensional solutions of the Yang-Baxter equation may 
be obtained.

\end{abstract}

\vfill
\vfill

\newpage
\pagestyle{plain}

The enveloping Lie algebra ${\cal U}(sl(2))$ has two distinct 
quantizations: The first one is called the Drinfeld-Jimbo 
deformation (standard $q$-deformation) \cite{Drin,Jim}, whereas the second 
one is called the Jordanian deformation (nonstandard $h$-deformation) 
\cite{Manin,Ohn} and may be obtained as a contraction of the 
Drinfeld-Jimbo ones \cite{Iran}. Recently there is much interest in studies 
relating to various aspects of the $h$-deformed algebra ${\cal U}_{h}(sl(2))$.
In particular, a two parametric deformation of the dual algebra 
${\cal F}un_{h,h'}(GL(2))$ 
was obtained in \cite{Agha}. This author also constructed \cite{Agha} 
the differential 
calculus in the quantum plane. Quantum de Rham complexes associated with 
the $h$-deformed algebra ${\cal F}un_{h}(sl(2))$ was given in 
\cite{Kiri}. The universal ${\cal R}$-matrix of the algebra 
${\cal U}_{h}(sl(2))$ was obtained 
\cite{Vlad,Shar,BH}. Various non-semisimple $h$-deformed algebras were 
constructed at contraction limits \cite{Ball,Shar,para}. The $h$-deformation 
was also extended to the case of supergroups \cite{Dab}.   

In \cite{Ohn} the fundamental representation of ${\cal U}_h(sl(2))$ was
obtained. The fundamental representation remains undeformed. 
The finite dimensional highest weight representations of ${\cal
U}_h(sl(2))$ were given in \cite{Dob} first by a direct
construction (cf. Proposition 3) and then by a
factor-representation of the corresponding Verma modules. The
latter was carried out by the standard singular vector
construction method, the new feature being the necessity for a
special construction of a homogeneous basis. In the direct
construction of \cite{Dob}, however, the representation action of
one of the generators was not given explicitly (except for the
three-dimensional irrep).
The main purpose of the present Letter is to
construct a nonlinear invertible map between the generators of 
${\cal U}_h(sl(2))$ and the classical $sl(2)$ generators. This immediately 
yields the irreducible representations of ${\cal U}_h(sl(2))$ in an explicit 
and simple manner.  The known universal ${\cal R}$-matrix \cite{Vlad,Shar,BH} 
of ${\cal U}_h(sl(2))$ algebra and our construction of its arbitrary 
irreducible representations may now be combined to generate new finite 
dimensional solutions of the Yang-Baxter equation.  

\smallskip
\smallskip  
   
Let $h$ be an arbitrary complex parameter. The algebra ${\cal U}_{h}(sl(2))$ 
is an associative algebra over $\CC$ generated by $H$, $T$, $T^{-1}$ and 
$Y$, satisfying the relations \cite{Ohn}
\begin{eqnarray}
  && TT^{-1}=T^{-1}T=1,   \\
  && [H,T] = T^{2}-1\;,\qquad\qquad\qquad\qquad [H,T^{-1}] = T^{-2}-1\;,  \\
  && [H,Y]=-{1\over 2}(YT+TY +YT^{-1}+T^{-1}Y), \\
  && [Y, T] = - {h\over 2}(HT+T H)\;, \qquad\;\qquad
  [Y, T^{-1}] = {h\over 2}(HT^{-1}+T^{-1}H)\;.  
\end{eqnarray}
We will not use the (non-standard) coalgebra structure 
\cite{Ohn} of ${\cal U}_h(sl(2))$ in our work.
Following \cite{Ohn}, we define a new generator $X:=h^{-1}\log{T}$. The 
algebra ${\cal U}_{h}(sl(2))$ is then an associative algebra over $\CC$
generated by $H$, $X$ and $Y$, satisfying the commutations relations 
\cite{Ohn}      
\begin{eqnarray}
  && [H,X]= 2 {\sinh hX \over h}, \\
  && [H,Y]=-Y(\cosh hX)-(\cosh hX)Y,\\  
  && [X,Y]=H. 
\end{eqnarray}
The Casimir element of ${\cal U}_{h}(sl(2))$ is given by \cite{Ball}  
\begin{eqnarray}
&&C={1\over 4h}\biggl(YT-YT^{-1}+TY-T^{-1}Y\biggr)+{1\over 4}H^2+
{1\over 16}(T^2+T^{-2}-2), \nonumber \\
&&\phantom{C}={1\over 2h}\biggl(Y(\sinh hX)+(\sinh hX)Y\biggr)+{1\over 4}H^2+
{1\over 4}(\sinh hX)^2.
\end{eqnarray}
The elements $Y^{k} H^{l} T^{s}$, $(k,l,s) \in \NN\times \NN\times \ZZ$ 
(resp. $Y^{k} H^{l} X^{s}$, $(k,l,s) \in \NN\times \NN\times \NN$)
build a Poincar\'e-Birkhoff-Witt basis of the algebra ${\cal U}_{h}(sl(2))$ 
\cite{Ohn}.  

The Verma module $M$ of ${\cal U}_{h}(sl(2))$ is generated by the highest 
weight vector $w_{0}$ with the highest weight $\lambda$         
\begin{eqnarray}
&& H.w_{\lambda}=\lambda\; w_{\lambda}, \qquad\qquad\qquad\qquad   
T.w_{\lambda}=\varepsilon\; w_{\lambda},\;\;\;\;(\varepsilon =\pm 1), \\
&& {1\over 2}(T-T^{-1})w_{\lambda}=0.
\end{eqnarray}
Let us start with the case $\varepsilon =1$. (The case $\varepsilon =-1$ will 
be included at the end.) Then the vectors $w_{m}:=Y^{m}.w_{0}$, $m>0$ and 
$w_{0}:=w_{\lambda}$ are a basis 
of $M$. If $\lambda =1$, then $w_{s}^{1/2}:=w_{2}$ 
is a primitive (singular) vector in $M$ ($X.w_{s}^{1/2}=0$). It is the 
first primitive vector in $M$ below $w_{\lambda}$. Therefore 
$L^{1\over 2}:=M/({\cal U}_{h}(sl(2)).w_{s}^{1/2})$ is a finite two 
($j={1\over 2}$) dimensional irreducible representation of 
${\cal U}_{h}(sl(2))$ and is given by \cite{Ohn}     
\begin{eqnarray}
H=\pmatrix{ 1 &  0\cr
		0& -1 \cr}, \qquad\qquad
X=\pmatrix{ 0 & 1 \cr
		0& 0 \cr},  \qquad\qquad
Y=\pmatrix{ 0& 0 \cr
		1& 0 \cr}. 
\end{eqnarray}
Similarly, if $\lambda=2$, then $w_{s}^{1}:=w_{3}+h^2 w_{1}$ is 
a singular vector in $M$. The quotient space $L^{1}:=M/({\cal U}_{h}
(sl(2)).w_{s}^{1})$ is a finite three ($j=1$) dimensional irreducible 
representation of ${\cal U}_{h}(sl(2))$ spanned by \cite{Dob}
\begin{eqnarray}
H=\pmatrix{2 &  0 & -2  h^2\cr
	0 & 0 & 0  \cr
		0 & 0 & -2  \cr}, \qquad
X=\pmatrix{ 0 & 2  & 0 \cr
	0 & 0 & 2 \cr
		0& 0 & 0 \cr},\qquad
Y=\pmatrix{ 0& 0  & 0\cr
	   1 & 0 & - h^2 \cr 
		0 &  1& 0 \cr}. 
\end{eqnarray}  
With an eye to the future development, we present here the matrix 
representation, 
where $H$ is diagonal 
\begin{eqnarray}
H=\pmatrix{2 &  0 & 0\cr
	0 & 0 & 0  \cr
		0 & 0 & -2  \cr}, \qquad
X=\pmatrix{ 0 & 2  & 0 \cr
	0 & 0 & 2 \cr
		0& 0 & 0 \cr}, \qquad
Y=\pmatrix{ 0& -\displaystyle {h^2\over 2}  & 0\cr
	   1 & 0 & - \displaystyle {h^2\over 2} \cr 
		0 &  1& 0 \cr}.
\end{eqnarray} 
For $\lambda=3$, the first singular vector in $M$ below $w_{\lambda}$ reads 
$w_{s}^{3/2}:=w_{4}+6h^2w_{2}$. Consequently, $L^{3\over 2}:=M/({\cal U}_{h}
(sl(2))$. $w_{s}^{3/2})$ is a finite four ($j={3\over 2}$) dimensional 
irreducible representation of ${\cal U}_{h}(sl(2))$ given by 
\begin{eqnarray}
H=\pmatrix{3 &  0 & -6  h^2 & 0\cr
	0 & 1 & 0 & -18h^2  \cr
	0 & 0 & -1 & 0  \cr
	0 & 0 & 0 & -3  \cr},\;\;\;
X=\pmatrix{ 0 & 3  & 0 & -6h^2 \cr
         0 & 0  & 4 & 0 \cr
	0 & 0 & 0 & 3 \cr
	 0& 0& 0 & 0 \cr},\;\;\;
Y=\pmatrix{ 0& 0 & 0 & 0\cr
	   1 & 0 & 0 & 0  \cr 
	0 & 1 & 0& -6h^2 \cr
	0 & 0 & 1& 0 \cr}. \;\;\;
\end{eqnarray} 
When $H$ is diagoalized, the representation $(14)$ assumes the form  
\begin{eqnarray}
H=\pmatrix{3 &  0 & 0 & 0\cr
	0 & 1 & 0 & 0 \cr
	0 & 0 & -1 & 0  \cr
	0 & 0 & 0 & -3  \cr},\;\;
X=\pmatrix{ 0 & 3  & 0 & 3h^2 \cr
         0 & 0  & 4 & 0 \cr
	0 & 0 & 0 & 3 \cr
	 0& 0& 0 & 0 \cr},\;\;
Y=\pmatrix{ 0& -\displaystyle{3h^2\over 2} & 0 & \displaystyle{9h^4\over 4} \cr
	   1 & 0 & -3h^2 & 0  \cr 
	0 & 1 & 0& -\displaystyle{3h^2\over 2} \cr
	0 & 0 & 1& 0 \cr}.
\end{eqnarray}
The other finite dimensional irreducible representations associated to 
$\lambda=4,\;5,\;6,\cdots$ may be obtained in a similar way. The singular 
vectors for $\lambda \leq 3$ were obtained analyticaly 
in \cite{Dob}. For $\lambda\leq 6$, the singular vectors may be constructed 
from results derived in \cite{Dob} using REDUCE. Nonetheless, 
the computations involved for obtaining the explicit expressions at higher 
dimensions become quite complicated very rapidly. However, a glance at the 
above representations $(13)$ and $(15)$, where the operator $H$ has been 
diagonalized, seems to suggest a close kinship between the ${\cal U}_{q}
(sl(2))$ generators and their analogues for the classical $sl(2)$ algebra. 
We pursue this route here.      

\smallskip
\smallskip

Let $\lbrace w_{-j}^{j},\;w_{-j+1}^{j},\cdots,\;w_{m}^{j},\cdots,\;
w_{j-1}^{j},\;w_{j}^{j}\rbrace$ be a set of basis vectors for the 
irreducible quotient module $L_{j}:=M/({\cal U}_{h}(sl(2)).w_{s}^{j})$ 
of dimension $(2j+1)$, where $H$ is diagonal. We wish to point out that,
with respect to the preceding examples, we introduce 
in the following a change of convention leading, for $h=0$, to the standard 
$sl(2)$ representations,  where the usual hermiticity properties of the 
generators are satisfied. The action of the generators 
$H$ and $X$ on these basis vectors are given by 
\begin{eqnarray}
&& H.w_{m}^{j}=2\;m\;w_{m}^{j}, \\
&& X.w_{m}^{j}=\sum_{k=0}^{[(j-m-1)/2]}h^{2k}\;a^{m+2k+1}_{m}\;w_{m+2k+1}^{j},
\end{eqnarray}
where $[x]$ stands for the integer part of $x$. Using the commutation relation 
$(5)$, we can easily prove that 
\begin{eqnarray}   
&& a_{m}^{m+2k+1}=f_{k}\;a_{m}\;a_{m+1}\cdots a_{m+2k-1}\;a_{m+2k} 
\qquad\qquad k\geq 0,
\end{eqnarray}
where $a_{m}=\sqrt{(j-m)(j+m+1)}$ and the coefficients $f_k$, 
$k\in \NN$ obey the following recurrence relation
\begin{eqnarray}   
\displaystyle f_k=\left\{ \matrix{1 & \qquad  k=0 \cr
 & \cr
\displaystyle {1\over 2k}\displaystyle\sum_{s=1}^{k} {1\over (2s+1)!} 
\sum_{i_1+i_2+\cdots + i_{2s+1}=k-s}
f_{i_1}f_{i_2}\cdots f_{i_{2s+1}} & \qquad  k\geq 1. \cr}\right. 
\end{eqnarray}
To solve $(19)$, we introduce a generating function
defined as $F(x)=\sum_{k=0}^{\infty}f_{k} x^{2k}$. The relation (19) 
requires $F(x)$ to satisfy the following differential 
equation  
\begin{eqnarray}   
{d (xF(x))\over dx}={\sinh(xF(x))\over x}, \qquad\qquad F(0)=1, 
\end{eqnarray} 
whose solution is given by 
\begin{eqnarray}
\displaystyle F(x)={2\over x}\;\hbox{arctanh}
({x\over 2})=\sum_{k=0}^{\infty}{1\over (2k+1)\;2^{2k}}\;x^{2k}.
\end{eqnarray} 
The solution for $\lbrace f_{k}\;|\; k\geq 0\rbrace$ now reads
\begin{eqnarray}   
f_k= {1\over (2k+1)2^{2k}}.
\end{eqnarray} 
Using the relations $(18)$ and $(22)$, the action of the generator $X$ on
the basis may be written as  
\begin{eqnarray}   
X.w_{m}^{j}=\sum_{k=0}^{[(j-m-1)/2]}{h^{2k}\over 
(2k+1)\;2^{2k}}\;\biggl(\;\prod_{s=0}^{2k}\;a_{m+s}\;\biggr)\; 
w_{m+2k+1}^{j}
\end{eqnarray}
and, evidently
\begin{eqnarray}
T.w_{m}^{j}= \displaystyle\left\lbrace\matrix{
w_{m}^{j} + \displaystyle \sum_{k\geq 1}^{j-m}{h^{k}\over 2^{k-1}}\; 
 \Biggl(\displaystyle \prod_{s=0}^{k-1}a_{m+s}\Biggr) \; w_{m+k}^{j}
   &\hbox{if}& \hfill m<j \cr
& & \cr
&& \cr
w_{j}^{j}&\hbox{if}& \hfill m=j \cr}\right.
\end{eqnarray}

Now we demonstrate the previously mentioned invertible map between the 
generators ${\cal U}_h(sl(2))$ algebra and the generators of the classical 
$sl(2)$. Let $J_{+}$, $J_{-}$ and $J_{3}$ be the classical $sl(2)$ generators 
acting on the basis $w_{m}^{j}$, $-j\leq m\leq j$ as 
\begin{eqnarray}
J_{+}.w_{m}^{j}=a_{m}\;w_{m+1}^{j},\qquad J_{-}.w_{m}^{j}=
a_{m-1}\;w_{m-1}^{j}, \qquad J_{3}.w_{m}^{j}=m\;w_{m}^{j}. 
\end{eqnarray}
Then the previous analysis indicates that the actions of the generators 
$X$ and $H$ on the vector space are equivalent to the following classical 
constructs   
\begin{eqnarray}   
X = {2\over h}\;\hbox{arctanh}({h\;J_{+}\over 2}),
\qquad\qquad\qquad\qquad H=2\;J_{3}
\end{eqnarray}
and
\begin{eqnarray}
\displaystyle T={1\;+\;\displaystyle {h\;J_{+}\over 2}\over 1\;-\;
\displaystyle {h\;J_{+}\over 2}}= 
\;1 \;+\;\displaystyle \sum_{k\geq 1}\;{h^{k}\over 2^{k-1}}\;J_{+}^{k}. 
\end{eqnarray}  
The inverse of the map $(25)$ and $(27)$ reads
\begin{eqnarray}   
J_{+} = {2\over h}\;\tanh({h\;X\over 2})={2\over h} 
\;\Biggl({T\;-\;1\over T\;+\;1}\Biggr),
\qquad\qquad\qquad J_{3}={H\over 2}.
\end{eqnarray}
To obtain a similar map for the generator $Y$, we take an ansatz   
\begin{eqnarray}   
Y=\varphi_{h}(J_{+})\;J_{-}\;\varphi_{h}(J_{+})  
\end{eqnarray}
with the natural condition $\varphi_{0}(x)=1$. The commutation 
relation $(6)$ and the useful identity $[J_{3},\; \phi(J_{+})]=
\phi'(J_{+})\;J_{+}$ now require  
the function $\varphi_{h}(x)$ to satisfy the following differential 
equation
\begin{eqnarray}
{1\over \varphi_{h}(x)}{\partial \varphi_{h}(x)\over \partial x} = 
-{h^2 x/4\over 1-h^2x^2/4}, 
\end{eqnarray}
whose solution is given by $\varphi_{h}(x)=(1-h^2x^2/4)^{1/2}$. The 
generator $Y$ may be written finally as  
\begin{eqnarray}   
&& Y=\sqrt{1-{h^2J_{+}^2 \over 4}}\;J_{-}\; \sqrt{1-{h^2J_{+}^2 \over 4}}, \\
&& \phantom{Y}=\sum_{k=0}^{\infty}{(-1)^{k}\;h^{2k} \over 
2^{2k}}\;\biggl(\sum_{s=0}^{k}\zeta_{s}\;\zeta_{k-s}
\;J_{+}^{2s}\;J_{-}
\;J_{+}^{2k-2s}\biggr),
\end{eqnarray}  
where, $\zeta_{0}=1$ and $\zeta_{k}= \displaystyle
{1\over k!}\prod_{s=0}^{k-1}({1\over 2}-s)$ for $k\geq 1$. The inverse map 
of $(31)$ readily follows 
\begin{eqnarray}
&& J_{-}=\cosh({h\;X\over 2})\;Y\;\cosh({h\;X\over 2}), \\
&& \phantom{ J_{-}}= {1\over 4}\;(T^{1/2}+T^{-1/2})\;Y \;(T^{1/2}+T^{-1/2}).
\end{eqnarray} 
The generator $Y$ act on the basis $w_{m}^{j}$ as  
\begin{eqnarray}   
&& Y.w_{m}^{j}= \sum_{k=0}^{[(j-m+1)/2]}{(-1)^{k}\;h^{2k} \over 
2^{2k}}\;\biggl(\sum_{s=0}^{k}\zeta_{s}\;\zeta_{k-s}
\;\psi_{k-s}{(m+2s-1)}\;a_{m+2s-1}\;\psi_{2s}(m)\biggr) w_{m+2k-1}^{j},
\nonumber \\
&& 
\end{eqnarray}  
where $\psi_{0}(m)=1$ and $\psi_{s}(m)=\prod_{k=0}^{2s-1}a_{m+k}$ for 
$s\geq 1$. The expressions $(26)$, $(27)$ and $(31)$ constitute the 
realization of the Jordanian algebra ${\cal U}_{h}(sl(2))$ with the 
classical generators via a nonlinear map. This map gives the representation 
characterized by $\varepsilon =1$. It may be directly verified that the 
commutation relations $(5)$, $(6)$ and $(7)$ are satisfied. Expressed in 
terms of the classical $sl(2)$ generators, the Casimir element $(8)$ of 
${\cal U}_{h}(sl(2))$ is just equal to 
\begin{eqnarray}
C={1\over 2}\;(J_{+}J_{-}+J_{-}J_{+})+J_{3}^2
\end{eqnarray}
and its value is equal to the classical one
\begin{eqnarray}
C.w_{m}^{j}=j(j+1)\;w_{m}^{j},\qquad\qquad -j \leq m \leq j. 
\end{eqnarray}    
Obviously the representation $(16)$, $(23)$, $(24)$ and $(35)$ 
($\varepsilon =1$) is a 
$h$-analog to the spin $j$ representation of $sl(2)$.

Let us mention that there is a $\CC$-algebra automorphism 
$\omega$ of ${\cal U}_{h}(sl(2))$ such that 
\begin{eqnarray}  
&& \omega(T)=T^{-1},\qquad \qquad \qquad \qquad 
\qquad\omega (T^{-1})=T,\nonumber \\
&&  \omega (Y)=-Y, \qquad \qquad \qquad \qquad \qquad 
\omega(H)=H
\end{eqnarray}
and, evidently  
\begin{eqnarray}  
\omega (X)=-X.
\end{eqnarray}
(For $h=0$, this automorphism reduces to the classical one 
$(J_{+},J_{-},J_{3})\longrightarrow (-J_{+},-J_{-},J_{3})$). 
Also there is a second $\CC$-algebra automorphism $\varpi$ of 
${\cal U}_{h}(sl(2))$ defined as
\begin{eqnarray}  
&& \varpi (T)=-\;T,\qquad \qquad \qquad \qquad 
\qquad\varpi (T^{-1})=-\;T^{-1},\nonumber \\
&&  \varpi (Y)=-Y, \qquad \qquad \qquad \qquad \qquad 
\varpi(H)=-H
\end{eqnarray}
and, evidently  
\begin{eqnarray}  
\varpi(X)=X+ {i\;\pi \over h}.
\end{eqnarray}
The representation induced by the automorphism $\varpi$ is characterized by 
$X.w_{j}^{j}={i\;\pi \over h}\;w_{j}^{j}$ (instead of $X.w_{j}^{j}=0$). The 
representations characterized by $\varepsilon =-1$ in $(9)$ are simply 
obtained from those presented in this paper using the automorphism $\varpi$.
These representations have evidently no classical ($h=0$) limit.  
(See \cite{Rosso} for the corresponding role of $\varepsilon$ for 
the ${\cal U}_q(sl(2))$ case.) 
Finally, the map associated to the irreducible representations associated 
to $\varepsilon$ ($\varepsilon =\pm 1$) is described by    
\begin{eqnarray}
&& T=  \varepsilon\;{1\;+\;\displaystyle {h\;J_{+} \over 2}\over  
1\;-\;\displaystyle{h\;J_{+} \over 2}}, \\   
&& Y = \varepsilon\;\sqrt{1-{h^2\;J_{+}^{2}\over 4}}\;J_{-}\;
\sqrt{1-{h^2\;J_{+}^{2}\over 4}},\\
&& H=2\;\varepsilon\;J_{3}
\end{eqnarray} 
and, where
\begin{eqnarray}
X={(1-\varepsilon )\;i\;\pi\;\over 2\;h}+{2
\over h}\;\hbox{arctanh}({h\;J_{+}\over 2}).
\end{eqnarray} 
This implies for the Casimir element $(8)$ ($\varepsilon =\pm 1$)
\begin{eqnarray}
&&C={1\over 4h}\biggl(YT-YT^{-1}+TY-T^{-1}Y\biggr)+{1\over 4}H^2+
{1\over 16}(T^2+T^{-2}-2), \nonumber \\
&&\phantom{C}={1\over 2h}\biggl(Y(\sinh hX)+(\sinh hX)Y\biggr)+{1\over 4}H^2+
{1\over 4}(\sinh hX)^2, \nonumber \\
&& \phantom{C}={1\over 2}(J_{+}\;J_{-}+J_{-}\;J_{+})+J_{3}^{2}.
\end{eqnarray}

The main result achieved in this Letter is the realization of the ${\cal 
U}_h(sl(2))$ generators in terms of the undeformed $sl(2)$ generators. This
facilitates immediate construction of the representation of the ${\cal 
U}_h(sl(2))$ algebra. It is interesting to note that, via this map, the 
generators $(H,X,Y)$ satisfying the algebraic structure $(5)$, $(6)$ and
$(7)$, may also be equipped with an induced cocommutative coproduct. 
This way, the generators $(H, X,Y)$ are viewed as composite elements of the 
${\cal U}(sl(2))$ algebra.  Similarly, the undeformed generators $(J_\pm , 
J_3)$, via the inverse map, may be viewed as elements of the 
${\cal U}_h(sl(2))$ algebra; and, thus, may be endowed with an induced 
noncocommutative coproduct \cite{abdess}. An example of nonlinear mapping
between $sl(2)$ and ${\cal U}_{q}(sl(2))$ can be found in \cite{Zachos}. 
This map is also discussed in \cite{abdess}. Applications of our formalism 
to a type of deformation of the $3$-dimensional euclidean space has been 
presented in \cite{ACC}.

Using our previous construction of the arbitrary irreducible representations 
of the ${\cal U}_h(sl(2))$ algebra via the map, we now obtain the finite 
dimensional solutions of the Yang-Baxter equation.  To this end, we use the 
expression of the universal ${\cal R}$-matrix of the ${\cal U}_h(sl(2))$ 
algebra, derived in \cite{BH}: 
\begin{eqnarray}
{\cal R} = \exp \left\{ -h X \otimes e^{hX} H \right\} 
           \exp \left\{ h e^{hX} H \otimes X \right\}\,. 
\end{eqnarray}
Our construction of the map (42)-(45) suggests that an arbitrary finite 
dimensional irreducible representation of the ${\cal U}_h(sl(2))$ algebra 
may be characterized by the indices $(j,\varepsilon )$, where $(2j+1) \in \ZZ _+$ 
and $\varepsilon = \pm 1$.  The finite dimensional representations of the 
$R$-matrix now satisfy the Yang-Baxter equation 
\begin{eqnarray}
R^{j_1,\varepsilon_1;j_2,\varepsilon_2}_{12} 
R^{j_1,\varepsilon_1;j_3,\varepsilon_3}_{13} 
R^{j_2,\varepsilon_2;j_3,\varepsilon_3}_{23} = 
R^{j_2,\varepsilon_2;j_3,\varepsilon_3}_{23} 
R^{j_1,\varepsilon_1;j_3,\varepsilon_3}_{13}
R^{j_1,\varepsilon_1;j_2,\varepsilon_2}_{12}\,.
\end{eqnarray}
The first few solutions for the classical case $(\varepsilon = 1)$ read: 
\begin{eqnarray}
R^{\frac{1}{2},\varepsilon_1=1;1,\varepsilon_2=1} & = & \left( 
\begin{array}{cccccc}
1 & 2h & 2h^2 & -2h & 2h^2 & 0 \\
0 & 1 & 2h & 0 & 0 & 2h^2 \\
0 & 0 & 1 &  0 & 0 & 2h \\
0 & 0 & 0 & 1 & -2h & 2h^2 \\
0 & 0 & 0 & 0 & 1 & -2h \\
0 & 0 & 0 & 0 & 0 & 1 
\end{array} \right)\,, \\
R^{\frac{1}{2},\varepsilon_1=1;\frac{3}{2},\varepsilon_2=1} & = & \left( 
\begin{array}{cccccccc}
1 & 3h & 6h^2 & 9h^3 & -3h & 3h^2 & 0 & 9h^4 \\
0 & 1 & 4h & 6h^2 & 0 & -h & 4h^2 & 0 \\
0 & 0 & 1 & 3h & 0 & 0 & h & 3h^2 \\
0 & 0 & 0 & 1 & 0 & 0 & 0 & 3h \\
0 & 0 & 0 & 0 & 1 & -3h & 6h^2 & -9h^3 \\
0 & 0 & 0 & 0 & 0 & 1 & - 4h & 6h^2 \\
0 & 0 & 0 & 0 & 0 & 0 & 1 & -3h \\
0 & 0 & 0 & 0 & 0 & 0 & 0 & 1 
\end{array} \right)\,.
\end{eqnarray}
Of special interest are the `non-classical' solutions of the type 
\begin{eqnarray}
R^{\frac{1}{2},\varepsilon_1 = 1;\frac{1}{2}, \varepsilon_2 = -1} = \left( 
\begin{array}{cccc}
-1 & h & h & h^2 \\
0 & -1 & 0 & -h \\
0 & 0 & -1 & -h \\
0 & 0 & 0 & -1 
\end{array} \right)\,.  
\end{eqnarray}
All such finite dimensional solutions of the Yang-Baxter equation may be 
similarly obtained from our construction of the irreducible representations of 
the ${\cal U}_h(sl(2))$ algebra and its known universal ${\cal R}$-matrix 
\cite{BH}. 

\vskip 2cm 
\noindent {\bf Acknowledgments:} 
 
We thank Claude de Calan and Jean Lascoux for interesting discussions. We 
also thank Ramaswamy Jagannathan for pointing out certain references to us. 
One of us (RC) wants to thank A. Chakrabarti for a kind invitation. He is 
also grateful to the members of the CPTH group for their kind hospitality.

\end{document}